\documentclass[12pt]{article}
\usepackage[T1]{fontenc}
\usepackage[utf8]{inputenc}
\usepackage{lmodern,amsmath,amssymb,amsthm}
\usepackage[a4paper,margin=1in]{geometry}
\usepackage{microtype,indentfirst}
\usepackage{xcolor}
\usepackage[colorlinks=true,linkcolor=blue,citecolor=blue,urlcolor=blue]{hyperref}
\hypersetup{pdftitle={Sharp singular Moser-Trudinger inequalities and n-Laplacian equations with degenerate potentials}}
\numberwithin{equation}{section}
\newtheorem{theorem}{Theorem}[section]
\newtheorem{lemma}[theorem]{Lemma}
\newtheorem{proposition}[theorem]{Proposition}
\newtheorem{corollary}[theorem]{Corollary}
\theoremstyle{remark}
\newtheorem{remark}[theorem]{Remark}
\newcommand{\R}{\mathbb R}
\newcommand{\dd}{\,dx}
\newcommand{\doi}[1]{\href{https://doi.org/#1}{\nolinkurl{doi: #1}}}

\title{Sharp singular Moser--Trudinger inequalities and\linebreak
$n$-Laplacian equations with degenerate potentials}
\author{%
Hao Hu\thanks{School of Mathematical Sciences, Chongqing Normal University,
Chongqing 401331, China.\newline
Email: \href{mailto:huh0201@126.com}{huh0201@126.com}.}
\and
Yanjun Liu\thanks{School of Mathematical Sciences, Chongqing Normal University,
Chongqing 401331, China.\newline
Email: \href{mailto:liuyj@mail.nankai.edu.cn}{liuyj@mail.nankai.edu.cn}.}
\and
Huabo Zhang\thanks{School of Mathematics and Statistics, Nantong University,
Nantong 226019, China.\newline
Email: \href{mailto:lutzhb@163.com}{lutzhb@163.com}.}
}
\date{}

\begin{document}
\maketitle

\begin{abstract}
Let $n\geq2$, $0<\beta<n$, and let $V$ be a continuous nonnegative potential on $\R^n$ with a
positive lower bound outside a ball. Under the constraint $\int_{\R^n}(|\nabla
u|^n+V|u|^n)\dd\leq1$, we prove the singular Moser--Trudinger inequality with sharp coefficient
$(1-\beta/n)\alpha_n$, where $\alpha_n=n|\mathbb S^{n-1}|^{1/(n-1)}$. The subcritical supremum is
attained for every such potential. If $V$ is radial and $r^\beta V(r)$ is nondecreasing, the
critical supremum is also attained. This includes potentials that vanish on a ball and decrease to a
positive limit at infinity. Under the same radial hypothesis, a bound for the tangential gradient
shows that every extremal is radial and that the tangential gradients of maximizing sequences
converge to zero in $L^n$, including at the critical coefficient. For the singular $n$-Laplacian
equation with critical exponential growth, we obtain an explicit sufficient level for Palais--Smale
compactness. Under the Ambrosetti--Rabinowitz condition and a smallness assumption at zero, this
gives a nonzero nonnegative mountain-pass solution for all sufficiently large coupling parameters.
\end{abstract}

{\small
\noindent\textbf{Keywords.} Singular Moser--Trudinger inequality, degenerate potential, extremal function, tangential gradient, $n$-Laplacian.\par
\smallskip
\noindent\textbf{2020 MSC.} 46E35, 35J92, 35B33, 35A15.
}

\section{Introduction and main results}
\label{sec:intro}

The Moser--Trudinger inequality describes the limiting Sobolev embedding when the gradient exponent
equals the dimension. On a bounded domain $\Omega\subset\R^n$, functions in $W_0^{1,n}(\Omega)$ need
not be bounded, but they remain exponentially integrable. Building on Trudinger's
work~\cite{Trudinger}, Moser~\cite{Moser} identified the optimal coefficient in
\begin{equation}
\sup_{\substack{u\in W_0^{1,n}(\Omega)\\ \|\nabla u\|_n\leq1}}
\int_\Omega e^{\alpha |u|^{n/(n-1)}}\dd<\infty
\quad\Longleftrightarrow\quad \alpha\leq\alpha_n,
\qquad
\alpha_n=n\omega_{n-1}^{1/(n-1)},
\label{eq:boundedMT}
\end{equation}
where $\omega_{n-1}$ denotes the surface measure of the unit sphere. Establishing the existence of
extremal functions requires additional analysis because compactness fails at the endpoint. The work
of Carleson--Chang~\cite{CarlesonChang} and Lin~\cite{Lin}, together with Lions's
concentration--compactness principle~\cite{Lions}, provides fundamental tools for this analysis.

On the whole space, the infinite volume necessitates subtracting the lower-order terms in the
exponential expansion. Throughout the paper, we write
\begin{equation}
 \Phi_n(t)=e^t-\sum_{j=0}^{n-2}\frac{t^j}{j!},\quad t\geq0.
\label{eq:Phi}
\end{equation}
The whole-space theory was developed by Cao~\cite{Cao}, Adachi--Tanaka~\cite{AdachiTanaka}, and
do \'O~\cite{doO}. Ruf~\cite{Ruf} and Li--Ruf~\cite{LiRuf} established sharp endpoint
inequalities under a full Sobolev norm constraint, while Ishiwata~\cite{Ishiwata} examined how the
constraint and the possible loss of compactness at infinity affect the existence of maximizers.

Singular weights introduce a further concentration mechanism. Adimurthi and
Sandeep~\cite{AdimurthiSandeep} proved that the weight $|x|^{-\beta}$ changes the bounded-domain
critical coefficient to
\begin{equation}
 \alpha_{n,\beta}=\left(1-\frac{\beta}{n}\right)\alpha_n,
 \qquad 0<\beta<n.
\label{eq:sharpconstant}
\end{equation}
Adimurthi--Yang~\cite{AdimurthiYang} established the corresponding whole-space inequality under a
full Sobolev norm and applied it to nonhomogeneous $n$-Laplacian equations. For the constant
potential $V\equiv\tau>0$, Li--Yang~\cite{LiYangSingular} proved attainment throughout the range
$0<\alpha\leq\alpha_{n,\beta}$. Nguyen--Takahashi~\cite[Theorems~1 and~2]{NguyenTakahashi}
studied the full norm with potential $|x|^{-\gamma}$. In their radial theorems, $\gamma=-\sigma<0$
gives $V(x)=|x|^\sigma$. For this potential the energy decreases under Schwarz rearrangement, while
the integral with weight $|x|^{-\beta}$ increases. Their radial result therefore also gives the
supremum over the full space for this power potential.

On bounded domains, Csat\'o--Nguyen--Roy~\cite{CsatoRoyNguyen} obtained singular extremals by
$n$-harmonic transplantation. Anisotropic versions replace the Euclidean gradient norm by a convex
homogeneous norm. Wang--Xia~\cite{WangXia} studied blow-up for a planar Finsler--Liouville
equation. Zhou--Zhou~\cite{ZhouZhouBounded,ZhouZhouWhole} treated anisotropic Dirichlet and full
Sobolev constraints. For singular or exact-growth inequalities on the whole space we refer to
Liu~\cite{LiuFinsler,LiuExactGrowth}, Li~\cite{LiAnisotropic} and Guo--Liu~\cite{GuoLiuWhole}.
Refinements on bounded domains appear in~\cite{LiuImproved,GuoLiuBounded}, and equations with the
Finsler $n$-Laplacian and critical exponential growth are studied in~\cite{DasTiwari}.

When $V$ vanishes on an open set, the potential term alone does not control the $L^n$ norm. Chen--Lu--Zhu~\cite{ChenLuZhu,ChenLuZhu2} obtained sharp unweighted inequalities for such potentials. In
dimension two, Zhang--Zhu~\cite[Theorem~1.1]{ZhangZhu} already proved the singular endpoint
estimate when the potential vanishes on a ball and is uniformly positive outside a larger ball. A
lower comparison potential gives the two-dimensional upper bound of the present paper directly from
their theorem. In every dimension the upper bound also follows by combining the power change of
variables used in~\cite{NguyenTakahashi} with the unweighted estimate
in~\cite[Theorem~1.1]{ChenLuZhu2}. Sun--Song--Zou~\cite{SunSongZou} proved unweighted endpoint
attainment for radial nondecreasing degenerate potentials, with an additional $L^p$ term in the
exponent whose coefficient may be zero. Related refinements were obtained in~\cite{LuZhu,Zhu}. Xue--Xu--Zhu~\cite{XueXuZhu} studied singular extremals on $\R^n$ in the presence of a trapping
potential.

Our question is which radial potentials allow the unweighted extremal to be transferred to the
singular problem. Write $s=r^{1-\beta/n}$. The potential in the transformed problem is a positive
constant times $r^\beta V(r)$. Its monotonicity suffices to rearrange functions in the new
variables, even when $V$ itself decreases on an interval. Keeping the angular derivatives in the
change of variables gives an estimate for the tangential gradient of every competitor in terms of
its distance in functional value from the supremum. At the critical coefficient this estimate still
applies to every maximizing sequence.

For the equation, Yang~\cite{YangEquations} obtained positive weak solutions on $\R^n$, and Lam--Lu~\cite{LamLuEquations} proved existence and multiplicity results for nonuniformly elliptic
equations with singular nonlinearities. Zhang--Chen~\cite[Theorems~1.2 and~1.4]{ZhangChen} proved
compactness in $L^q(\R^n,|x|^{-\beta}\dd)$ for $q\geq n$ and obtained positive ground states for
uniformly positive potentials under their hypotheses on the nonlinearity. Zhang--Zhu~\cite[Theorem~1.4]{ZhangZhu} treated a planar degenerate potential under additional trapping and
growth conditions. For unweighted equations, Chen--Qin--R\u{a}dulescu--Tang~\cite{ChenQinRadulescuTang} obtained ground states and mountain-pass solutions for several
classes of potentials. Here the weight gives the compactness needed for a mountain-pass solution at
large coupling. The potential may vanish on a compact set and may oscillate at infinity.

Throughout the paper, $n\geq2$, $0<\beta<n$, and $n'=n/(n-1)$. We assume that
\begin{equation*}
 V\in C(\R^n,[0,\infty)),\qquad
 \inf_{|x|\geq R_0} V(x)\geq v_0
 \quad\text{for some }R_0,v_0>0.
\tag{$V_1$}\label{ass:V}
\end{equation*}
We use the energy space
\begin{equation}
 E=\left\{u\in W^{1,n}(\R^n):\int_{\R^n}V(x)|u|^n\dd<\infty\right\},
 \qquad
 \|u\|^n=\int_{\R^n}\bigl(|\nabla u|^n+V(x)|u|^n\bigr)\dd.
\label{eq:space}
\end{equation}
We write $E_W$ and $\|\cdot\|_W$ for the same space and norm with $V$ replaced by a potential $W$.
In particular, $E=E_V$ and $\|\cdot\|=\|\cdot\|_V$. Lemma~\ref{lem:coercivity} shows that this norm
controls the $W^{1,n}$ norm. For $\alpha>0$, set
\begin{equation}
 I_\alpha(u)=\int_{\R^n}|x|^{-\beta}\Phi_n(\alpha|u|^{n'})\dd,\qquad
 S_\alpha=\sup_{\substack{u\in E\\\|u\|\leq1}} I_\alpha(u).
\label{eq:sup}
\end{equation}

\begin{theorem}\label{thm:sharp}
Assume \eqref{ass:V} and $0<\beta<n$. Then
\begin{equation}
 S_\alpha<\infty\quad\text{for }0<\alpha\leq\alpha_{n,\beta},
 \qquad
 S_\alpha=\infty\quad\text{for }\alpha>\alpha_{n,\beta}.
\label{eq:sharp}
\end{equation}
Moreover, for every fixed $u\in E$ and every $\alpha>0$,
\begin{equation}
 \int_{\R^n}\frac{\Phi_n(\alpha|u|^{n'})}{|x|^\beta}\dd<\infty.
\label{eq:individual}
\end{equation}
The sharp coefficient does not depend on the size of the bounded region in which $V$ may vanish.
\end{theorem}

\begin{remark}\label{rem:normalization}
The critical coefficient must be interpreted together with the energy constraint. For any $\rho>0$,
the substitution $u=\rho v$ gives
\begin{equation}
 \sup_{\|u\|\leq\rho}
 \int_{\R^n}|x|^{-\beta}\Phi_n(\alpha|u|^{n'})\dd
 =S_{\alpha\rho^{n'}}.
\label{eq:radiusnormalization}
\end{equation}
Consequently, this supremum is finite exactly when
$\alpha\rho^{n'}\leq\alpha_{n,\beta}$. The coefficient in Theorem~\ref{thm:sharp} is independent of the particular potential, although the finite value of the supremum may depend on $V$. An estimate obtained solely from an equivalence of norms can give a smaller admissible coefficient.
\end{remark}

The hypothesis used for critical attainment is
\begin{equation*}
 V(x)=V(r),\quad r=|x|,\qquad
 r\longmapsto r^\beta V(r)\quad\text{is nondecreasing on }[0,\infty).
\tag{$V_{2,\beta}$}\label{ass:radial}
\end{equation*}
Set $a=1-\beta/n$. For $x\neq0$ and $\theta=x/|x|$, write
$\partial_r u(x)=\nabla u(x)\cdot\theta$ and
$\nabla_\tau u(x)=\nabla u(x)-\partial_r u(x)\theta$ for the radial derivative and the tangential component of the gradient.

\begin{theorem}\label{thm:extremal}
Assume \eqref{ass:V}. For each fixed $0<\alpha<\alpha_{n,\beta}$, $S_\alpha$ is attained, and every
maximizing sequence in the unit ball has a subsequence converging strongly in $E$ to a maximizer.

If \eqref{ass:radial} also holds, $S_{\alpha_{n,\beta}}$ is attained. For each
$0<\alpha\leq\alpha_{n,\beta}$ and every $u\in E$ with $\|u\|\leq1$,
\begin{equation}
 \frac{I_\alpha(u)}{S_\alpha}
 +(1-a^n)\int_{\R^n}|\nabla_\tau u|^n\dd
 \leq \|u\|^n .
 \label{eq:tangentialbound}
\end{equation}
Hence every maximizer under \eqref{ass:radial} is radial about the origin. A maximizer may be chosen
nonnegative and nonincreasing in $r$. For any maximizing sequence $(u_k)$ in the unit ball at such a
fixed $\alpha$, $\nabla_\tau u_k\to0$ in $L^n(\R^n)$. All maximizers considered above have norm one.

A nonnegative maximizer $u_\alpha$ satisfies, in the weak sense,
\begin{equation}
 -\Delta_n u_\alpha+V(x)u_\alpha^{n-1}
 =\frac{|x|^{-\beta}u_\alpha^{1/(n-1)}
 \Phi_n'(\alpha u_\alpha^{n'})}
 {\displaystyle\int_{\R^n}|x|^{-\beta}u_\alpha^{n'}
 \Phi_n'(\alpha u_\alpha^{n'})\dd},
\label{eq:extremalEL}
\end{equation}
where $\Delta_n u=\operatorname{div}(|\nabla u|^{n-2}\nabla u)$.
\end{theorem}

\begin{remark}\label{rem:criticalscope}
For a potential that decreases outside a ball, let
$\chi(r)=\min\{\max\{r-1,0\},1\}$ and set
\begin{equation}
 V(r)=\chi(r)\left(1+\frac{1}{1+r^\beta}\right).
 \label{eq:decreasingexample}
\end{equation}
This potential vanishes on $B_1$ and decreases strictly to $1$ for $r\geq2$. Nevertheless,
\[
 r^\beta V(r)=\chi(r)\left(r^\beta+\frac{r^\beta}{1+r^\beta}\right)
\]
is the product of two nonnegative nondecreasing functions. Hence \eqref{ass:V} and
\eqref{ass:radial} hold, and the critical supremum for \eqref{eq:decreasingexample} is attained.

For a singularity at $x_0\in\R^n$, the change of variables $y=x-x_0$ preserves \eqref{ass:V} after
the exterior ball is enlarged. The sharp bound, subcritical attainment and the equation result
therefore hold for the weight $|x-x_0|^{-\beta}$. The critical attainment and symmetry conclusions
apply if $V(x)=U(|x-x_0|)$ with $r^\beta U(r)$ nondecreasing.
\end{remark}

The exponential estimate also yields an existence result for the singular equation
\begin{equation}
 -\Delta_n u+V(x)|u|^{n-2}u
 =\lambda\frac{f(u)}{|x|^\beta}\quad\text{in }\R^n,
 \qquad u\in E,
\label{eq:problem}
\end{equation}
where $\lambda>0$. Write $F(t)=\int_0^t f(s)\,ds$ and assume that $f$ satisfies the following conditions.
\begin{align*}
 &f\in C(\R),\quad f(t)=0\text{ for }t\leq0,\quad
 \lim_{t\downarrow0}\frac{f(t)}{t^{n-1}}=0.
 \tag{$f_1$}\label{ass:f1}\\
 &\text{there exists }\alpha_0>0\text{ such that}\quad
 \lim_{t\to\infty}\frac{f(t)}{e^{b t^{n'}}}
 =\begin{cases}0,&b>\alpha_0,\\+\infty,&0<b<\alpha_0\end{cases}.
 \tag{$f_2$}\label{ass:f2}\\
 &\text{there exists }\mu>n\text{ such that}\quad
 0<\mu F(t)\leq t f(t)\quad\text{for all }t>0.
 \tag{$f_3$}\label{ass:f3}
\end{align*}
The parameter $\alpha_0$ specifies the critical growth of $f$ and does not bound the energy norm.
The smallness condition required for compactness will be imposed on the mountain-pass level.

\begin{theorem}\label{thm:equation}
Assume \eqref{ass:V}, $0<\beta<n$, and \eqref{ass:f1}--\eqref{ass:f3}. There exists $\lambda_0>0$
such that, for every $\lambda>\lambda_0$, problem~\eqref{eq:problem} has a nonzero nonnegative weak
solution. It is a mountain-pass critical point of
\begin{equation}
 J_\lambda(u)=\frac{1}{n}\|u\|^n
 -\lambda\int_{\R^n}\frac{F(u)}{|x|^\beta}\dd.
\label{eq:energy}
\end{equation}
A sufficient condition for $J_\lambda$ to satisfy the Palais--Smale condition at level $c$ is
\begin{equation}
 0<c<c_*:=\left(\frac1n-\frac1\mu\right)
 \left(\frac{\alpha_{n,\beta}}{\alpha_0}\right)^{n-1}.
\label{eq:PSthreshold}
\end{equation}
The mountain-pass level $c_\lambda$ tends to zero as $\lambda\to\infty$.
\end{theorem}

\begin{corollary}\label{cor:model}
Let $p>n, 0<\beta<n$ and $\alpha_0>0$. Under \eqref{ass:V}, the equation
\begin{equation}
 -\Delta_n u+V(x)u^{n-1}
 =\lambda |x|^{-\beta}u^{p-1}e^{\alpha_0u^{n'}},
 \qquad u\geq0,\quad u\in E,
\label{eq:model}
\end{equation}
has a nonzero weak solution for all sufficiently large $\lambda$. Fix a nonzero $\varphi\in
C_c^\infty(\R^n)$ with $\varphi\geq0$. Choosing a sufficiently large multiple of $\varphi$ as the
endpoint in the mountain-pass construction gives
\begin{equation}
 c_\lambda\leq
 \left(\frac1n-\frac1p\right)
 \frac{\|\varphi\|^{np/(p-n)}}
 {\left(\lambda\displaystyle\int_{\R^n}|x|^{-\beta}\varphi^p\dd\right)^{n/(p-n)}}.
\label{eq:modelbound}
\end{equation}
\end{corollary}

Section~\ref{sec:prelim} gives the energy and compactness estimates. Section~\ref{sec:sharp}
computes the effect of the radial change of variables on an arbitrary function and proves the sharp
inequality. Section~\ref{sec:extremal} treats attainment and symmetry. Section~\ref{sec:equation}
proves the Palais--Smale estimate and the large-parameter existence result.

\section{Energy estimates and weighted compactness}
\label{sec:prelim}

\begin{lemma}\label{lem:coercivity}
Under \eqref{ass:V}, $E$ is a uniformly convex Banach space and
\begin{equation}
 \|u\|_q\leq C_q\|u\|\quad\text{for every }n\leq q<\infty.
\label{eq:Lqcontrol}
\end{equation}
Moreover, $C_c^\infty(\R^n)$ is dense in $E$, and bounded sequences in $E$ have subsequences
converging strongly in $L^q_{\mathrm{loc}}(\R^n)$ for every finite $q\geq1$.
\end{lemma}

\begin{proof}
Choose $R\geq R_0$ and a cutoff $\eta\in C_c^\infty(B_{2R})$ such that $0\leq\eta\leq1$, $\eta=1$ on
$B_R$, and $|\nabla\eta|\leq C/R$. Applying the Poincar\'e inequality to $\eta u\in
W_0^{1,n}(B_{2R})$ gives
\begin{align}
 \int_{B_R}|u|^n\dd
 &\leq C_R\int_{B_{2R}}|\nabla(\eta u)|^n\dd\notag\\
 &\leq C_R\int_{B_{2R}}|\nabla u|^n\dd
   +C_R\int_{B_{2R}\setminus B_R}|u|^n\dd\notag\\
 &\leq C_R\int_{\R^n}|\nabla u|^n\dd
   +\frac{C_R}{v_0}\int_{\R^n}V(x)|u|^n\dd.
\label{eq:cutoffLn}
\end{align}
The last inequality holds because the derivative of the cutoff is supported where $V\geq v_0$. Outside $B_R$,
\begin{equation}
 \int_{\R^n\setminus B_R}|u|^n\dd
 \leq\frac1{v_0}\int_{\R^n}V(x)|u|^n\dd.
\label{eq:exteriorLn}
\end{equation}
Combining these estimates yields $\|u\|_n\leq C\|u\|$. The whole-space Sobolev embeddings
$W^{1,n}(\R^n)\hookrightarrow L^q(\R^n)$ for $n\leq q<\infty$ then yield \eqref{eq:Lqcontrol}.

If $(u_k)$ is Cauchy in $E$, it converges in $W^{1,n}$, and $(V^{1/n}u_k)$ converges in $L^n$. An
almost-everywhere convergent subsequence shows that the latter limit is $V^{1/n}u$, proving
completeness. The map $u\mapsto(\nabla u,V^{1/n}u)$ embeds $E$ isometrically as a closed subspace of
a uniformly convex $L^n$ product space. Hence $E$ is uniformly convex. The Rellich theorem on
bounded balls gives the asserted local compactness.

To prove density, choose a cutoff satisfying $\eta_R=1$ on $B_R$, supported in $B_{2R}$, with
$0\leq\eta_R\leq1$ and $|\nabla\eta_R|\leq C/R$. The product rule gives
\[
 \|u-\eta_Ru\|^n
 \leq C\int_{|x|>R}(|\nabla u|^n+V|u|^n)\dd
       +CR^{-n}\int_{|x|>R}|u|^n\dd\longrightarrow0.
\]
For fixed $R$, mollification of $\eta_Ru$ gives convergence in $W^{1,n}$. If the mollifier radii are
less than one, all approximations are supported in $B_{2R+1}$, where $V$ is bounded. The potential
error is therefore bounded by $\sup_{B_{2R+1}}V$ times the $n$th power of the $L^n$ error. First
taking $R$ sufficiently large and then the mollifier radius sufficiently small proves density in
$E$. Radial cutoffs and mollifiers preserve radial symmetry, so the same argument proves density of
radial smooth compactly supported functions in the radial subspace.
\end{proof}

Lemma~\ref{lem:coercivity}, together with the decay of the singular weight, yields compact weighted
embeddings. Following Zhang--Chen~\cite{ZhangChen}, we include the proof to clarify the role of
the potential.

\begin{lemma}\label{lem:weighted}
For every $n\leq q<\infty$, the embedding
\begin{equation}
 E\hookrightarrow L^q(\R^n,|x|^{-\beta}\dd)
\label{eq:weightedembed}
\end{equation}
is continuous and compact. In particular, $u_k\rightharpoonup u$ in $E$ implies, after extraction if necessary,
\begin{equation}
 \int_{\R^n}|x|^{-\beta}|u_k-u|^q\dd\longrightarrow0.
\label{eq:weightedstrong}
\end{equation}
\end{lemma}

\begin{proof}
Choose $r>1$ such that $\beta r<n$ and set $r'=r/(r-1)$. H\"older's inequality and
Lemma~\ref{lem:coercivity} yield
\begin{equation}
 \int_{B_1}|x|^{-\beta}|u|^q\dd
 \leq\left(\int_{B_1}|x|^{-\beta r}\dd\right)^{1/r}
       \|u\|_{qr'}^q
 \leq C\|u\|^q.
\label{eq:weightedlocal}
\end{equation}
On $\R^n\setminus B_1$, the weight is at most one, so \eqref{eq:Lqcontrol} bounds the remaining
integral and establishes continuity.

To prove compactness, use local strong convergence in $L^{qr'}$ and H\"older's inequality to obtain
weighted convergence on each fixed ball. For $R\geq1$, the complementary integral satisfies
\begin{equation}
 \int_{|x|>R}|x|^{-\beta}|u_k-u|^q\dd
 \leq R^{-\beta}\|u_k-u\|_q^q\leq C R^{-\beta}.
\label{eq:polynomialtail}
\end{equation}
Letting $k\to\infty$ and then $R\to\infty$ proves \eqref{eq:weightedstrong}.
\end{proof}

\begin{remark}\label{rem:betarange}
The bound $\beta<n$ ensures local integrability, and $\beta>0$ gives decay at infinity. If
$\beta\geq n$, a smooth compactly supported function that is positive and constant near the origin
has an infinite weighted exponential integral for every $\alpha>0$, even after multiplication by a
small positive constant to satisfy the energy constraint. Thus the upper bound is necessary even for
smooth test functions.

To see why the lower bound is needed for compactness, take $V\equiv1$, a nonzero $\varphi\in
C_c^\infty(\R^n)$, and points $x_k$ with $|x_k|\to\infty$. Then $u_k(x)=\varphi(x-x_k)$ is bounded
in $E$ and converges weakly to zero, whereas $\|u_k\|_q=\|\varphi\|_q$. For $\beta>0$ and
sufficiently large $k$,
\begin{equation}
 \int_{\R^n}|x|^{-\beta}|u_k|^q\dd
 \leq C|x_k|^{-\beta}\|\varphi\|_q^q\longrightarrow0.
\label{eq:translatedbump}
\end{equation}
The singular weight therefore restores compactness in all weighted $L^q$ spaces covered by
Lemma~\ref{lem:weighted}, including $q=n$, without the stronger assumption $V(x)\to\infty$.
Compactness at critical exponential growth requires a separate argument.
\end{remark}

We shall use the following unweighted attainment theorem.

\begin{proposition}\label{prop:SSZ}
Suppose that $W\in C(\R^n,[0,\infty))$ is radial and nondecreasing, and
\begin{equation}
 0\leq\inf_{\R^n}W<\lim_{|x|\to\infty}W(x)\leq\infty.
\label{eq:SSZhyp}
\end{equation}
Then
\begin{equation}
 \sup_{\int_{\R^n}(|\nabla v|^n+W|v|^n)\dd\leq1}
 \int_{\R^n}\Phi_n(\alpha_n|v|^{n'})\dd<\infty,
\label{eq:SSZbound}
\end{equation}
and the supremum is attained. A maximizer can be chosen to be nonnegative, radial and nonincreasing,
with energy one.
\end{proposition}

Proposition~\ref{prop:SSZ} follows by setting the coefficient of the additional norm term equal to
zero in~\cite[Theorem~1.1]{SunSongZou}. The theorem includes the case
$\lim_{|x|\to\infty}W(x)=\infty$.

We record the rearrangement inequalities used below. If $W$ is radial and nondecreasing, and $w^*$
denotes the Schwarz rearrangement of $|w|$, then
\begin{equation}
 \int_{\R^n}|\nabla w^*|^n\dd\leq\int_{\R^n}|\nabla w|^n\dd,\qquad
 \int_{\R^n} W|w^*|^n\dd\leq\int_{\R^n}W|w|^n\dd.
 \label{eq:rearrangeenergy}
\end{equation}
The first inequality is the P\'olya--Szeg\H{o} inequality. For the second, among sets of a
prescribed measure a centered ball minimizes the integral of $W$. Integration over the level sets of
$|w|^n$ gives the claim, also for unbounded $W$ by monotone convergence. Equimeasurability gives,
for $b>0$,
\begin{equation}
 \int_{\R^n}\Phi_n(b|w^*|^{n'})\dd
 =\int_{\R^n}\Phi_n(b|w|^{n'})\dd .
 \label{eq:rearrangeobjective}
\end{equation}
For these rearrangement facts see~\cite{Kesavan}.

\begin{lemma}\label{lem:unweighted}
Under \eqref{ass:V},
\begin{equation}
 \sup_{\|u\|\leq1}\int_{\R^n}\Phi_n(\alpha_n|u|^{n'})\dd<\infty.
\label{eq:unweighted}
\end{equation}
Consequently, for $0<\alpha\leq\alpha_{n,\beta}$ and $R\geq1$,
\begin{equation}
 \sup_{\|u\|\leq1}\int_{|x|>R}
 |x|^{-\beta}\Phi_n(\alpha|u|^{n'})\dd\leq C R^{-\beta}.
\label{eq:exponentialtail}
\end{equation}
\end{lemma}

\begin{proof}
Take $W(x)=0$ for $|x|\leq R_0$ and $W(x)=v_0/2$ for $|x|\geq2R_0$, with linear interpolation in
$|x|$ on the intervening annulus. Then $0\leq W\leq V$ and $\|u\|_W\leq\|u\|$. The unweighted bound
for $W$ in~\cite[Theorem~1.1]{ChenLuZhu2} gives \eqref{eq:unweighted}. For $|x|>R$ we have
$|x|^{-\beta}\leq R^{-\beta}$. Since $\alpha\leq\alpha_{n,\beta}<\alpha_n$,
\eqref{eq:exponentialtail} follows from the monotonicity of $\Phi_n$.
\end{proof}

\section{The sharp singular inequality}
\label{sec:sharp}

We use the power change of radius from Nguyen--Takahashi~\cite[Section~2]{NguyenTakahashi}, up to
a constant dilation. Keeping the angular variable gives the following identities for a general
potential.

\begin{lemma}\label{lem:transform}
Suppose that \eqref{ass:V} holds. For $r>0$, $\theta\in\mathbb S^{n-1}$ and $a=1-\beta/n$, set
\begin{equation}
 s=r^a,\qquad v(s\theta)=a^{1/n'}u(r\theta),\qquad
 \widetilde V(s\theta)=a^{-n}r^\beta V(r\theta).
 \label{eq:transform}
\end{equation}
Put $\widetilde V(0)=0$. Then $\widetilde V$ is continuous and tends to infinity at infinity. The
map $u\mapsto v$ is a bounded linear bijection from $E$ to $E_{\widetilde V}$ with bounded inverse.
It satisfies
\begin{align}
 \int_{\R^n}|\nabla v(y)|^n\,dy
 &=\int_{\R^n}\bigl(|\partial_r u|^2+a^2|\nabla_\tau u|^2\bigr)^{n/2}\dd ,
 \label{eq:gradientidentity}\\
 \int_{\R^n}\widetilde V(y)|v(y)|^n\,dy
 &=\int_{\R^n}V(x)|u(x)|^n\dd ,
 \label{eq:potentialidentity}\\
 I_\alpha(u)&=\frac1a\int_{\R^n}
       \Phi_n\left(\frac{\alpha}{a}|v(y)|^{n'}\right)\,dy .
 \label{eq:objectiveidentity}
\end{align}
The equality in \eqref{eq:objectiveidentity} holds also for infinite integrals. The difference of
the two energies satisfies
\begin{equation}
 \|u\|^n-\|v\|_{\widetilde V}^n
 \geq(1-a^n)\int_{\R^n}|\nabla_\tau u|^n\dd .
 \label{eq:angulardifference}
\end{equation}
In particular, the energies agree for radial functions.
\end{lemma}

\begin{proof}
Write $x=r\theta$ and $y=s\theta$. The formula for $\widetilde V$ at $s>0$ can also be written as
$a^{-n}s^{n/a-n}V(s^{1/a}\theta)$. Since $n/a-n>0$ and $V$ is bounded near zero, this expression
tends uniformly to zero as $s\downarrow0$. For $s\geq R_0^a$ it is bounded below by
$a^{-n}v_0s^{n/a-n}$, which tends to infinity.

First let $u$ be smooth and compactly supported. For $r>0$, the Jacobian is
\[
 dy=a r^{n(a-1)}\dd=a r^{-\beta}\dd .
\]
Since $|v|^n=a^{n-1}|u|^n$ and $|v|^{n'}=a|u|^{n'}$, substitution gives \eqref{eq:potentialidentity}
and \eqref{eq:objectiveidentity}. Denote by $\nabla_{\mathbb S}$ the gradient on the unit sphere,
acting on the angular variable. On each ray the chain rule gives
\[
 \partial_s v=a^{-1/n}r^{1-a}\partial_r u,\qquad
 \frac{\nabla_{\mathbb S}v}{s}
      =a^{1/n'}r^{1-a}\frac{\nabla_{\mathbb S}u}{r}.
\]
Using $|\nabla v|^2=|\partial_s v|^2+s^{-2}|\nabla_{\mathbb S}v|^2$ and $|\nabla_\tau
u|=r^{-1}|\nabla_{\mathbb S}u|$ proves \eqref{eq:gradientidentity}. In particular,
\[
 a^n\int_{\R^n}|\nabla u|^n\dd\leq\int_{\R^n}|\nabla v|^n\,dy
 \leq\int_{\R^n}|\nabla u|^n\dd .
\]
For $B\geq0$, the function $A\mapsto(A+B)^{n/2}-(A+a^2B)^{n/2}$ is nondecreasing on $[0,\infty)$,
since $n/2\geq1$. Its value at $A=0$ is $(1-a^n)B^{n/2}$. Taking $A=|\partial_r u|^2$,
$B=|\nabla_\tau u|^2$ and integrating proves \eqref{eq:angulardifference}.

The forward map preserves membership in $L^n$. Indeed, the change of measure and Lemma~\ref{lem:weighted} give
\[
 \int_{\R^n}|v|^n\,dy=a^n\int_{\R^n}|x|^{-\beta}|u|^n\dd<\infty .
\]
For the inverse map take $R\geq\max\{1,R_0\}$. On $B_R$, use $s^{n/a-n}\leq R^\beta$ for $s\leq
R^a$. On its complement use $V\geq v_0$. The two estimates give
\[
 \int_{\R^n}|u|^n\dd
 \leq a^{-n}R^\beta\int_{|y|\leq R^a}|v|^n\,dy
       +v_0^{-1}\int_{\R^n}\widetilde V|v|^n\,dy .
\]
For a smooth compactly supported input in either variable, its transform is bounded near the origin
and smooth away from it. Integration by parts outside a ball of radius $\varepsilon$ identifies the
computed gradient as its weak gradient on $\R^n$. The boundary term is $O(\varepsilon^{n-1})$
against each smooth test function, and the gradient is in $L^n$ by the preceding bounds.

The estimates give $\|v\|_{\widetilde V}\leq\|u\|$ in the forward direction and $\|u\|\leq
a^{-1}\|v\|_{\widetilde V}$ in the inverse direction. Both spaces are complete and have dense
subspaces of smooth compactly supported functions by Lemma~\ref{lem:coercivity}. The maps therefore
extend continuously. To identify each extension, pass to an almost-everywhere convergent subsequence
of smooth approximations in both spaces. On every annulus the change of variables is a smooth
diffeomorphism, so the limits obey \eqref{eq:transform} there. This identifies the two extensions as
inverses almost everywhere. The gradient and potential identities now pass to the limit in norm. The
identity for $I_\alpha$ follows by measurable substitution, including when either side is infinite.
\end{proof}

\begin{proof}[Proof of Theorem~\ref{thm:sharp}]
Let $\|u\|\leq1$ and let $v$ be given by \eqref{eq:transform}. The potential $\widetilde V$
satisfies \eqref{ass:V} and $\|v\|_{\widetilde V}\leq1$. Since $\alpha_{n,\beta}/a=\alpha_n$,
\eqref{eq:objectiveidentity} and Lemma~\ref{lem:unweighted} applied to $\widetilde V$ give
\[
 I_{\alpha_{n,\beta}}(u)\leq\frac1a
 \sup_{\|w\|_{\widetilde V}\leq1}
 \int_{\R^n}\Phi_n(\alpha_n|w|^{n'})\,dy<\infty .
\]
This proves the upper bound for $\alpha=\alpha_{n,\beta}$. Monotonicity of $\Phi_n$ gives the
subcritical bounds.

To establish sharpness, let $0<\varepsilon<e^{-1}$, set $L=\log(1/\varepsilon)$, and define
\begin{equation}
 m_\varepsilon(x)=\omega_{n-1}^{-1/n}
 \begin{cases}
 L^{(n-1)/n},&|x|\leq\varepsilon,\\[2pt]
 L^{-1/n}\log(1/|x|),&\varepsilon<|x|<1,\\[2pt]
 0,&|x|\geq1.
 \end{cases}
\label{eq:Moser}
\end{equation}
Direct integration gives
\begin{equation}
 \int_{\R^n}|\nabla m_\varepsilon|^n\dd=1,
 \qquad \int_{\R^n}|m_\varepsilon|^n\dd=O(L^{-1}).
\label{eq:Moserenergy}
\end{equation}
Indeed, the annular gradient integral is $L^{-1}\int_\varepsilon^1 dr/r=1$. In the $L^n$ integral,
the inner ball contributes $O(\varepsilon^nL^{n-1})$, and the annulus contributes at most
\[
 \frac{C}{L}\int_0^1(\log(1/r))^n r^{n-1}\,dr=O(L^{-1}).
\]
Because $V$ is bounded on $\overline B_1$,
\begin{equation}
 d_\varepsilon:=\int_{\R^n}V|m_\varepsilon|^n\dd=O(L^{-1}),
 \qquad u_\varepsilon=\frac{m_\varepsilon}{(1+d_\varepsilon)^{1/n}}
 \quad\text{satisfies }\|u_\varepsilon\|=1.
\label{eq:normalizeMoser}
\end{equation}
On $B_\varepsilon$, the exponent has the uniform expansion
\begin{equation}
 \alpha|u_\varepsilon|^{n'}
 =\alpha\omega_{n-1}^{-1/(n-1)}L+O(1).
\label{eq:Mosercore}
\end{equation}
For sufficiently large $t$, $\Phi_n(t)\geq e^t/2$. Hence, for small $\varepsilon$,
\begin{align}
 \int_{\R^n}|x|^{-\beta}\Phi_n(\alpha|u_\varepsilon|^{n'})\dd
 &\geq C\varepsilon^{n-\beta}
 \exp\bigl(\alpha\omega_{n-1}^{-1/(n-1)}L\bigr)\notag\\
 &=C\exp\left[\left(\alpha\omega_{n-1}^{-1/(n-1)}-(n-\beta)\right)L\right].
\label{eq:Moserdivergence}
\end{align}
The right-hand side tends to infinity when $\alpha>\alpha_{n,\beta}$, proving sharpness. Since the
argument requires only local boundedness of $V$, the conclusion holds whether or not $V$ vanishes at
the origin.

It remains to establish integrability for each individual function. For $s,t\geq0$, monotonicity gives
\begin{equation}
 \Phi_n(s+t)\leq\Phi_n(2s)+\Phi_n(2t).
\label{eq:Phisplit}
\end{equation}
Fix $u\in E$ and $\alpha>0$. By density, choose $\psi\in C_c^\infty(\R^n)$ so close to $u$ in $E$ that
\[
 2^{n'}\alpha\|u-\psi\|^{n'}\leq\alpha_{n,\beta}.
\]
Using $|u|^{n'}\leq2^{1/(n-1)}(|\psi|^{n'}+|u-\psi|^{n'})$ in \eqref{eq:Phisplit}, we obtain
\begin{equation}
 \Phi_n(\alpha|u|^{n'})
 \leq\Phi_n(2^{n'}\alpha|\psi|^{n'})
   +\Phi_n(2^{n'}\alpha|u-\psi|^{n'}).
\label{eq:individualsplit}
\end{equation}
The first term is bounded and compactly supported, and is therefore integrable against the locally
integrable weight $|x|^{-\beta}$. After normalization, the endpoint estimate bounds the integral of
the second term. This proves \eqref{eq:individual}.
\end{proof}

We record two consequences for use below. For every $r\geq1$, there exists $C_r$ such that
\begin{equation}
 \Phi_n(t)^r\leq C_r\Phi_n(rt),\qquad t\geq0.
\label{eq:Phipower}
\end{equation}
For $0\leq t\leq1$ use $\Phi_n(t)\asymp t^{n-1}$, and for $t\geq1$ compare with $e^t$. For any fixed
$u\in E$ and $A>0$, there is also a $\delta>0$ such that
\begin{equation}
 \sup_{\|h\|\leq\delta}
 \int_{\R^n}|x|^{-\beta}\Phi_n(A|u+h|^{n'})\dd<\infty.
\label{eq:localuniform}
\end{equation}
To verify uniformity, choose $\psi\in C_c^\infty(\R^n)$ and then $\delta>0$ so that
\[
 2^{n'}A(\|u-\psi\|+\delta)^{n'}
 \leq\alpha_{n,\beta}.
\]
For $\|h\|\leq\delta$, apply the splitting argument to $u+h=\psi+(u-\psi+h)$. The integral of the
first term is a fixed finite quantity. Since $\|u-\psi+h\|\leq\|u-\psi\|+\delta$, normalization
bounds the integral of the second term by $S_{\alpha_{n,\beta}}$. This proves
\eqref{eq:localuniform} for every fixed $A>0$. The neighborhood depends on $u$ and $A$, so the
estimate does not extend uniformly to an arbitrary bounded subset of $E$.

\section{Extremal functions}
\label{sec:extremal}

We begin by establishing weak continuity of the exponential functional in the subcritical regime.

\begin{lemma}\label{lem:subcriticalcontinuity}
Let $0<\alpha<\alpha_{n,\beta}$. If $u_k\rightharpoonup u$ in $E$ and $\|u_k\|\leq1$, then
\begin{equation}
 \int_{\R^n}|x|^{-\beta}\Phi_n(\alpha|u_k|^{n'})\dd
 \longrightarrow
 \int_{\R^n}|x|^{-\beta}\Phi_n(\alpha|u|^{n'})\dd.
\label{eq:subcritconvergence}
\end{equation}
\end{lemma}

\begin{proof}
Choose $r>1$ with $r\alpha<\alpha_{n,\beta}$. By \eqref{eq:Phipower} and Theorem~\ref{thm:sharp},
\begin{equation}
 \sup_k\int_{\R^n}|x|^{-\beta}
 \Phi_n(\alpha|u_k|^{n'})^r\dd<\infty.
\label{eq:uniformintegrability}
\end{equation}
After passing to a subsequence, we have $u_k\to u$ almost everywhere. Since the measure
$|x|^{-\beta}\dd$ is finite on $B_R$, \eqref{eq:uniformintegrability} and Vitali's theorem yield
convergence of the integrals over this ball. Lemma~\ref{lem:unweighted} bounds the tails uniformly
by $CR^{-\beta}$, including the tail of the limit since $\|u\|\leq1$. Letting $R\to\infty$ proves
\eqref{eq:subcritconvergence}. The convergence holds for the full sequence because every subsequence
has a further subsequence with the same limit.
\end{proof}

\begin{proof}[Proof of subcritical attainment in Theorem~\ref{thm:extremal}]
Fix $0<\alpha<\alpha_{n,\beta}$ and let $(u_k)$ be a maximizing sequence in the unit ball.
Reflexivity provides a subsequence converging weakly to some $u\in E$, and
Lemma~\ref{lem:subcriticalcontinuity} yields
\[
 \int_{\R^n}|x|^{-\beta}\Phi_n(\alpha|u|^{n'})\dd=S_\alpha.
\]
Testing with a normalized, nonzero, compactly supported smooth function shows that $S_\alpha>0$, and
hence $u\neq0$. If $\|u\|<1$, strict monotonicity of $\Phi_n$ on $(0,\infty)$ implies that replacing
$u$ by $u/\|u\|$ increases the integral, contradicting the definition of $S_\alpha$. Therefore
$\|u\|=1$, and weak lower semicontinuity together with $\|u_k\|\leq1$ gives $\|u_k\|\to1$. Uniform
convexity of $E$ then yields strong convergence. Replacing $u$ by $|u|$ preserves the energy and
$I_\alpha(u)$, so a nonnegative maximizer exists.
\end{proof}

At $\alpha=\alpha_{n,\beta}$ the weak continuity in Lemma~\ref{lem:subcriticalcontinuity} fails. In
fact, the normalized Moser functions $u_\varepsilon$ in \eqref{eq:normalizeMoser} converge weakly to
zero in $E$. They are bounded in $E$, and \eqref{eq:Moserenergy} forces every weak limit to vanish.
At $\alpha=\alpha_{n,\beta}$ the coefficient of $L$ in \eqref{eq:Moserdivergence} is zero, so
$\liminf_{\varepsilon\downarrow0}I_{\alpha_{n,\beta}}(u_\varepsilon)>0$.

\begin{proof}[Proof of the radial conclusions of Theorem~\ref{thm:extremal}]
Assume \eqref{ass:radial}. Then $\widetilde V(s)=a^{-n}r^\beta V(r)$ is radial and nondecreasing in
$s$. For $0<\alpha\leq\alpha_{n,\beta}$ we claim that
\begin{equation}
 S_\alpha=\frac1a\sup_{\substack{w\in E_{\widetilde V}\\
                          \|w\|_{\widetilde V}\leq1}}
        \int_{\R^n}\Phi_n\left(\frac\alpha a |w|^{n'}\right)\,dy .
 \label{eq:supidentity}
\end{equation}
Indeed, the transform of each $u$ in the original unit ball lies in the transformed unit ball, and
\eqref{eq:objectiveidentity} gives the upper bound. For the reverse inequality, let $w$ belong to
the transformed unit ball and take the Schwarz rearrangement $w^*$ of $|w|$. Apply
\eqref{eq:rearrangeenergy} and \eqref{eq:rearrangeobjective} with $W=\widetilde V$. Then $w^*$ is
still admissible and has the same exponential integral as $w$. Its inverse transform
$u(r)=a^{-1/n'}w^*(r^a)$ is radial, so \eqref{eq:gradientidentity} and \eqref{eq:potentialidentity}
give $\|u\|=\|w^*\|_{\widetilde V}\leq1$. Applying \eqref{eq:objectiveidentity} to $u$ proves
\eqref{eq:supidentity}. In particular, the rearrangement is taken with respect to the transformed
potential, for which monotonicity follows from \eqref{ass:radial}.

At $\alpha=\alpha_{n,\beta}$ the coefficient in the transformed problem is $\alpha_n$. Since
$\widetilde V(0)=0$ and $\widetilde V(y)\to\infty$ as $|y|\to\infty$, Proposition~\ref{prop:SSZ}
supplies a nonnegative radial nonincreasing extremal for the right-hand side of
\eqref{eq:supidentity}. Its inverse transform attains $S_{\alpha_{n,\beta}}$ and is nonincreasing in
$r$. For $0<\alpha<\alpha_{n,\beta}$, take any maximizer already obtained in the first part of the
proof. Its transform attains the right-hand side of \eqref{eq:supidentity}. Rearranging there and
taking the inverse again yields a nonnegative radial nonincreasing maximizer. At every exponent
under consideration, strict increase of $I_\alpha(tu)$ in $t>0$ for nonzero $u$ gives norm one at a
maximizer.

To prove \eqref{eq:tangentialbound}, take $\|u\|\leq1$, let $v$ be its transform, and put
$t=\|v\|_{\widetilde V}$. If $t=0$, then $u=0$. For $0<t\leq1$, the series defining $\Phi_n$ starts
at degree $n-1$. Thus for $z\geq0$,
\[
 \Phi_n(t^{n'}z)
 =\sum_{j=n-1}^\infty\frac{t^{n'j}z^j}{j!}
 \leq t^n\Phi_n(z).
\]
Apply this inequality with $z=(\alpha/a)|v/t|^{n'}$. By \eqref{eq:objectiveidentity} and \eqref{eq:supidentity},
\[
 I_\alpha(u)\leq\frac{t^n}{a}\int_{\R^n}
       \Phi_n\left(\frac{\alpha}{a}|v/t|^{n'}\right)\,dy
 \leq t^n S_\alpha .
\]
Now combine $I_\alpha(u)/S_\alpha\leq\|v\|_{\widetilde V}^n$ with \eqref{eq:angulardifference}. This
proves \eqref{eq:tangentialbound}. For a maximizer, the inequality gives $\nabla_\tau u=0$ almost
everywhere. On every annulus, Sobolev slicing in polar coordinates gives $u(r,\cdot)\in
W^{1,n}(\mathbb S^{n-1})$ for almost every $r$. Its spherical gradient vanishes, so it is constant
almost everywhere on the connected sphere. Thus $u$ is radial. For a maximizing sequence in the unit
ball we obtain
\[
 (1-a^n)\int_{\R^n}|\nabla_\tau u_k|^n\dd
 \leq 1-\frac{I_\alpha(u_k)}{S_\alpha}\longrightarrow0 .
\]
This last estimate includes the critical exponent.
\end{proof}

\begin{proof}[Proof of the Euler--Lagrange equation in Theorem~\ref{thm:extremal}]
For fixed $\alpha>0$, $I_\alpha$ is continuously differentiable on $E$. Indeed, the derivative of
$t\mapsto\Phi_n(\alpha|t|^{n'})$ is $O(|t|^{n-1})$ near zero. For $|t|\geq1$ a larger exponential
coefficient absorbs its polynomial factor. The local uniform estimate \eqref{eq:localuniform} and
the weighted embeddings then give continuity of the derivative by the truncation argument in
Lemma~\ref{lem:functional}.

Let $u=u_\alpha\geq0$ be a maximizer of norm one. Since the derivative of the constraint $\|u\|^n=1$
takes the value $n$ when applied to $u$, the Lagrange multiplier rule gives $\theta\in\R$ such that,
for every $h\in E$,
\begin{align}
 \frac{\alpha n}{n-1}\int_{\R^n}|x|^{-\beta}u^{1/(n-1)}
 \Phi_n'(\alpha u^{n'})h\dd
 =\theta n\int_{\R^n}
 \bigl(|\nabla u|^{n-2}\nabla u\cdot\nabla h+Vu^{n-1}h\bigr)\dd.
\label{eq:multiplier}
\end{align}
Testing with $h=u$ and using $\|u\|=1$ gives
\[
 \theta=\frac{\alpha}{n-1}\int_{\R^n}|x|^{-\beta}
 u^{n'}\Phi_n'(\alpha u^{n'})\dd>0.
\]
To see that the integral is finite, bound the factor multiplying $|x|^{-\beta}$ by $Cu^n$ on
$\{u\leq1\}$ and by $C\Phi_n(bu^{n'})$ on $\{u>1\}$ for any fixed $b>\alpha$. Now use
Lemma~\ref{lem:weighted} and \eqref{eq:individual}. It is positive because $u\neq0$ and
$\Phi_n'(t)>0$ for $t>0$. Substituting this expression into \eqref{eq:multiplier} yields
\eqref{eq:extremalEL} with the stated denominator.
\end{proof}

\begin{remark}\label{rem:smallalpha}
As $\alpha\downarrow0$, the leading-order behavior of the extremal problem is governed by the
weighted $L^n$ embedding. More precisely,
\begin{equation}
 S_\alpha=\frac{\alpha^{n-1}}{(n-1)!}
 \max_{\|u\|\leq1}\int_{\R^n}|x|^{-\beta}|u|^n\dd
 +O(\alpha^n)\quad\text{as }\alpha\downarrow0.
\label{eq:smallalpha}
\end{equation}
Lemma~\ref{lem:weighted} and weak compactness of the energy ball show that the maximum on the
right-hand side is positive, finite and attained. The scaling argument forces every maximizer to
have norm one.

To establish the expansion, fix $0<\alpha_1<\alpha_{n,\beta}$. The power series of $\Phi_n$ gives,
for $0<\alpha\leq\alpha_1$ and $t\geq0$,
\begin{equation}
 0\leq\Phi_n(\alpha t)-\frac{\alpha^{n-1}t^{n-1}}{(n-1)!}
 \leq\left(\frac\alpha{\alpha_1}\right)^n\Phi_n(\alpha_1t).
\label{eq:smallalpharemainder}
\end{equation}
Applying this estimate with $t=|u|^{n'}$ and integrating gives a uniform remainder bound of
$(\alpha/\alpha_1)^nS_{\alpha_1}$ on the unit energy ball, proving \eqref{eq:smallalpha}.

The maximizers also converge. Let $\alpha_k\downarrow0$ and choose nonnegative maximizers $u_k$ for
$S_{\alpha_k}$. The uniform remainder estimate in \eqref{eq:smallalpha} shows that their weighted
$L^n$ integrals converge to the maximum in that equation. Along a subsequence, $u_k\rightharpoonup
u_0$ in $E$ and $u_k\to u_0$ in weighted $L^n$. Thus $u_0$ attains the positive maximum and has norm
one. Uniform convexity gives strong convergence in $E$, and the Lagrange multiplier rule yields the
weighted eigenvalue equation
\begin{equation}
 -\Delta_n u_0+V(x)u_0^{n-1}
 =\frac{|x|^{-\beta}u_0^{n-1}}
 {\displaystyle\int_{\R^n}|x|^{-\beta}u_0^n\dd},
 \qquad \|u_0\|=1,
\label{eq:limitweightedproblem}
\end{equation}
in the weak sense. This limiting polynomial problem remains valid when the potential is nonradial or
vanishes on a bounded region. The convergence is subsequential and does not require uniqueness of
the limiting maximizer.
\end{remark}

\section{Singular \texorpdfstring{$n$}{n}-Laplacian equations}
\label{sec:equation}

A weak solution of \eqref{eq:problem} is a function $u\in E$ such that
\begin{equation}
 \int_{\R^n}\bigl(|\nabla u|^{n-2}\nabla u\cdot\nabla h
   +V(x)|u|^{n-2}uh\bigr)\dd
 =\lambda\int_{\R^n}|x|^{-\beta}f(u)h\dd
 \quad\text{for every }h\in E.
\label{eq:weaksolution}
\end{equation}
The following estimates establish the finiteness of both sides and provide the variational
framework. The growth coefficient $\alpha_0$ in \eqref{ass:f2} is distinct from the admissible
exponential coefficient on the unit energy ball.

\begin{lemma}\label{lem:functional}
Assume \eqref{ass:V}, \eqref{ass:f1} and \eqref{ass:f2}. The functional $J_\lambda$ in
\eqref{eq:energy} belongs to $C^1(E,\R)$, with
\begin{equation}
 J_\lambda'(u)h=
 \int_{\R^n}\bigl(|\nabla u|^{n-2}\nabla u\cdot\nabla h
 +V|u|^{n-2}uh\bigr)\dd
 -\lambda\int_{\R^n}|x|^{-\beta}f(u)h\dd.
\label{eq:derivative}
\end{equation}
For every $b>\alpha_0$, $p>n$ and $\varepsilon>0$, there is $C_\varepsilon>0$ such that
\begin{align}
 |f(t)|&\leq\varepsilon|t|^{n-1}
 +C_\varepsilon|t|^{p-1}\Phi_n(b|t|^{n'}),
\label{eq:fgrowth}\\
 |F(t)|&\leq\frac\varepsilon n|t|^n
 +C_\varepsilon|t|^p\Phi_n(b|t|^{n'}),\qquad t\in\R.
\label{eq:Fgrowth}
\end{align}
\end{lemma}

\begin{proof}
Condition \eqref{ass:f1} yields the first term in \eqref{eq:fgrowth} for small positive $t$. For
large $t$, choosing a coefficient strictly between $\alpha_0$ and $b$ in \eqref{ass:f2} yields an
exponential upper bound that can be absorbed into $t^{p-1}\Phi_n(bt^{n'})$. Continuity controls the
remaining bounded interval. Integrating and using the monotonicity of $\Phi_n$ gives
\eqref{eq:Fgrowth}. Both estimates hold for $t\leq0$ since $f(t)=F(t)=0$.

For every $u\in E$, \eqref{eq:individual} ensures that the product of $\Phi_n(b|u|^{n'})$ with any
nonnegative power of $|u|$ is integrable with respect to $|x|^{-\beta}\dd$. On $|u|\leq1$, this
follows from $\Phi_n(b|u|^{n'})\leq C|u|^n$. On $|u|>1$ a larger exponential coefficient absorbs the
power. Consequently, $F(u)$ is integrable. The same argument gives $f(u)\in
L^{n'}(\R^n,|x|^{-\beta}\dd)$, because $|f(u)|^{n'}\leq C|u|^n$ for $|u|\leq1$ and has controlled
exponential growth for $|u|>1$.

To prove continuity of the derivative, suppose $u_k\to u$ in $E$ and fix $s>n/(n-1)$. Applying
\eqref{eq:localuniform} at $u$ with $A=2sb$ bounds the weighted exponential integrals of $u_k$
uniformly for large $k$. Then \eqref{eq:fgrowth} with $\varepsilon=1$, H\"older's inequality and
\eqref{eq:Phipower} yield
\begin{align*}
 \int_{\R^n}|x|^{-\beta}|f(u_k)|^s\dd
 &\leq C\int_{\R^n}|x|^{-\beta}|u_k|^{s(n-1)}\dd\\
 &\quad+C\left(\int_{\R^n}|x|^{-\beta}|u_k|^{2s(p-1)}\dd\right)^{1/2}\\
 &\qquad{}\times\left(\int_{\R^n}|x|^{-\beta}
 \Phi_n(2sb|u_k|^{n'})\dd\right)^{1/2}.
\end{align*}
Both polynomial exponents are at least $n$, so Lemma~\ref{lem:weighted} and the exponential bound
provide a uniform bound on the right-hand side. On any fixed ball, the weighted measure is finite
and $s>n/(n-1)$ ensures uniform integrability of $|f(u_k)|^{n'}$. Along a subsequence, $u_k\to u$
almost everywhere. Continuity of $f$ and Vitali's theorem then imply convergence of $f(u_k)$ to
$f(u)$ in weighted $L^{n'}$ on that ball.

To control the exterior region, split it into $|u_k|\leq1$ and $|u_k|>1$. On the first set, use
$|f(u_k)|^{n'}\leq C|u_k|^n$. On the second, for any fixed $r>1$ and its conjugate exponent $r'$,
\eqref{eq:fgrowth} and H\"older's inequality yield a bound by a constant multiple of
\begin{equation}
 \int_{|x|>R}|x|^{-\beta}|u_k|^n\dd
 +\left(\int_{|x|>R}|x|^{-\beta}
 |u_k|^{\frac{n(p-1)r'}{n-1}}\dd\right)^{1/r'}.
\label{eq:nonlinearitytail}
\end{equation}
The exponential factor is uniformly bounded by \eqref{eq:Phipower} and \eqref{eq:localuniform}.
Since $|x|^{-\beta}\leq R^{-\beta}$ outside $B_R$, the unweighted $L^q$ estimates in
Lemma~\ref{lem:coercivity} bound \eqref{eq:nonlinearitytail} by $C(R^{-\beta}+R^{-\beta/r'})$. The
same estimate holds for $u$, so the tails vanish uniformly as $R\to\infty$. Together with local
convergence, this proves global convergence in weighted $L^{n'}$, first along subsequences and then
along the full sequence. The continuous weighted $L^n$ embedding then gives continuity of the
nonlinear derivative in $E^*$.

The Fr\'echet remainder is estimated by
\begin{align*}
 &\left|\int_{\R^n}|x|^{-\beta}
   \bigl(F(u+h)-F(u)-f(u)h\bigr)\dd\right|\\
 &\quad\leq C\|h\|\sup_{0\leq t\leq1}
 \left(\int_{\R^n}|x|^{-\beta}
 |f(u+th)-f(u)|^{n'}\dd\right)^{(n-1)/n}
 =o(\|h\|).
\end{align*}
The inequality follows from the fundamental theorem of calculus and H\"older's inequality. The
supremum tends to zero by the continuity just proved and the bound $\|th\|\leq\|h\|$. The energy
term is continuously differentiable, being the sum of the $n$th powers of the $L^n$ norms of $\nabla
u$ and $V^{1/n}u$, proving \eqref{eq:derivative}. The same continuity argument applies to the
exponential functional in Section~\ref{sec:extremal}, since it requires only an $O(|t|^{n-1})$ bound
near zero and exponential upper bounds at infinity.
\end{proof}

We next establish mountain-pass geometry and estimate the associated minimax level.

\begin{lemma}\label{lem:geometry}
Under the hypotheses of Theorem~\ref{thm:equation}, for every $\lambda>0$ there exist
$\rho_\lambda,b_\lambda>0$ and $e_\lambda\in E$ such that
\begin{equation}
 J_\lambda(u)\geq b_\lambda\quad\text{if }\|u\|=\rho_\lambda,
 \qquad \|e_\lambda\|>\rho_\lambda,\quad J_\lambda(e_\lambda)<0.
\label{eq:geometry}
\end{equation}
For the corresponding mountain-pass construction, the endpoint may be chosen on the ray through a
fixed nonzero nonnegative $\varphi\in C_c^\infty(\R^n)$. The resulting levels satisfy $c_\lambda>0$
and $c_\lambda\to0$ as $\lambda\to\infty$.
\end{lemma}

\begin{proof}
Fix $b>\alpha_0$ and $p>n$. For $\|u\|\leq\rho$ with $2b\rho^{n'}\leq\alpha_{n,\beta}$,
\eqref{eq:Phipower}, Theorem~\ref{thm:sharp}, and the weighted embeddings give
\begin{align}
 \int_{\R^n}|x|^{-\beta}|u|^p\Phi_n(b|u|^{n'})\dd
 &\leq\left(\int_{\R^n}|x|^{-\beta}|u|^{2p}\dd\right)^{1/2}
 \left(\int_{\R^n}|x|^{-\beta}\Phi_n(b|u|^{n'})^2\dd\right)^{1/2}\notag\\
 &\leq C\|u\|^p.
\label{eq:smallnormestimate}
\end{align}
Consequently,
\[
 J_\lambda(u)\geq\left(\frac1n-C\lambda\varepsilon\right)\|u\|^n
 -C_\varepsilon\lambda\|u\|^p.
\]
Choose $\varepsilon$ so that $C\lambda\varepsilon\leq1/(2n)$ and then $0<\rho_\lambda\leq\rho$ such
that $C_\varepsilon\lambda\rho_\lambda^{p-n}\leq1/(4n)$. These choices are possible because $p>n$
and yield $J_\lambda(u)\geq\rho_\lambda^n/(4n)$ on $\|u\|=\rho_\lambda$, so we may take
$b_\lambda=\rho_\lambda^n/(4n)$.

Condition \eqref{ass:f3} implies that $F(t)/t^\mu$ is nondecreasing on $(0,\infty)$, because
\[
 \frac{d}{dt}\frac{F(t)}{t^\mu}
 =\frac{tf(t)-\mu F(t)}{t^{\mu+1}}\geq0.
\]
In particular, $F(t)\geq F(1)t^\mu$ for $t\geq1$. Fix a nonzero $\varphi\in C_c^\infty(\R^n)$ with
$\varphi\geq0$, and choose $\delta>0$ such that $A=\{x:\varphi(x)\geq\delta\}$ has positive measure.
Since $\beta<n$ and $A$ is bounded, $0<\int_A|x|^{-\beta}\dd<\infty$. For $t\geq1/\delta$,
\[
 \int_{\R^n}|x|^{-\beta}F(t\varphi)\dd
 \geq F(1)\delta^\mu t^\mu\int_A|x|^{-\beta}\dd,
\]
and consequently
\begin{equation}
 J_\lambda(t\varphi)\leq\frac{t^n}{n}\|\varphi\|^n-C\lambda t^\mu
 \longrightarrow-\infty\quad\text{as }t\to\infty.
\label{eq:negativepath}
\end{equation}
We may therefore choose an endpoint $e_\lambda$ on this ray and define the mountain-pass value by
\begin{equation}
 c_\lambda=\inf_{\gamma\in\Gamma_\lambda}\max_{t\in[0,1]}J_\lambda(\gamma(t)),
 \quad
 \Gamma_\lambda=\{\gamma\in C([0,1],E):\gamma(0)=0,\ \gamma(1)=e_\lambda\}.
\label{eq:mountainpass}
\end{equation}
Every such path crosses the sphere on which the functional is positive, so $c_\lambda\geq b_\lambda>0$.

For $\lambda\geq1$, choose a fixed, sufficiently large endpoint $T\varphi$ with $J_1(T\varphi)<0$,
and decrease $\rho_\lambda$ if necessary so that $\rho_\lambda<\|T\varphi\|$. Since $F\geq0$, this
is a valid endpoint for every larger value of $\lambda$. The straight path gives
\begin{equation}
 0<c_\lambda\leq\max_{t\geq0}J_\lambda(t\varphi).
\label{eq:raylevel}
\end{equation}
Continuity and \eqref{eq:negativepath} ensure the existence of a maximizing point $t_\lambda$. Since
$t=0$ is admissible, $J_1(t_\lambda\varphi)\geq J_\lambda(t_\lambda\varphi)\geq0$. The limit
$J_1(t\varphi)\to-\infty$ therefore confines all $t_\lambda$ to a bounded interval independent of
$\lambda\geq1$. If $t_\lambda\geq\delta>0$ along a sequence with $\lambda\to\infty$, then
\[
 J_\lambda(t_\lambda\varphi)
 \leq C-\lambda\int_{\R^n}|x|^{-\beta}F(\delta\varphi)\dd
 \longrightarrow-\infty,
\]
contradicting the nonnegativity of the maximum. Here \eqref{ass:f3} ensures that $F$ is increasing
and positive on $(0,\infty)$. Thus $t_\lambda\to0$ and
\[
 0<c_\lambda\leq J_\lambda(t_\lambda\varphi)
 \leq\frac{t_\lambda^n}{n}\|\varphi\|^n\longrightarrow0.
\]
This proves the level estimate.
\end{proof}

The norm bound supplied by \eqref{ass:f3} now yields compactness below the sufficient threshold in
\eqref{eq:PSthreshold}.

\begin{lemma}\label{lem:PS}
For every fixed $\lambda>0$, $J_\lambda$ satisfies the Palais--Smale condition at each level $0<c<c_*$.
\end{lemma}

\begin{proof}
Let $J_\lambda(u_k)\to c$ and $J_\lambda'(u_k)\to0$ in $E^*$. By \eqref{ass:f3},
\begin{align}
 \left(\frac1n-\frac1\mu\right)\|u_k\|^n
 &\leq J_\lambda(u_k)-\frac1\mu J_\lambda'(u_k)u_k\notag\\
 &=c+o(1)+o(1)\|u_k\|.
\label{eq:PSnorm}
\end{align}
If $\|u_k\|$ were unbounded, dividing \eqref{eq:PSnorm} by $\|u_k\|^n$ along an unbounded
subsequence would yield $1/n-1/\mu\leq0$, contradicting $\mu>n$. Thus $(u_k)$ is bounded, the term
$o(1)\|u_k\|$ tends to zero, and
\begin{equation}
 \limsup_{k\to\infty}\|u_k\|^n
 \leq\frac{c}{1/n-1/\mu}
 <\left(\frac{\alpha_{n,\beta}}{\alpha_0}\right)^{n-1}.
\label{eq:PSsmallnorm}
\end{equation}
The strict inequality allows us to choose
\[
 \limsup_k\|u_k\|<M<
 \left(\frac{\alpha_{n,\beta}}{\alpha_0}\right)^{(n-1)/n}.
\]
Choose $b>\alpha_0$ with $bM^{n'}<\alpha_{n,\beta}$ and then $r>1$ sufficiently close to one that
\begin{equation}
 rbM^{n'}<\alpha_{n,\beta}.
\label{eq:PSmargin}
\end{equation}
After discarding finitely many terms if necessary, Theorem~\ref{thm:sharp} applied to $u_k/\|u_k\|$,
together with \eqref{eq:Phipower}, yields
\begin{equation}
 \sup_k\int_{\R^n}|x|^{-\beta}\Phi_n(b|u_k|^{n'})^r\dd<\infty.
\label{eq:PSexp}
\end{equation}
Along a subsequence, Lemma~\ref{lem:weighted} gives $u_k\rightharpoonup u$ in $E$ and strong
convergence in every weighted $L^q$, $q\geq n$. To estimate the nonlinear term tested against
$u_k-u$, fix $p>n$ and set $r'=r/(r-1)$. By \eqref{eq:fgrowth}, \eqref{eq:PSexp} and H\"older's
inequality,
\begin{align}
 &\left|\int_{\R^n}|x|^{-\beta}f(u_k)(u_k-u)\dd\right|\notag\\
 &\quad\leq C
 \left(\int_{\R^n}|x|^{-\beta}|u_k|^n\dd\right)^{(n-1)/n}
 \left(\int_{\R^n}|x|^{-\beta}|u_k-u|^n\dd\right)^{1/n}\notag\\
 &\qquad+C\left(\int_{\R^n}|x|^{-\beta}|u_k|^{pr'}\dd\right)^{(p-1)/(pr')}
 \left(\int_{\R^n}|x|^{-\beta}|u_k-u|^{pr'}\dd\right)^{1/(pr')}
 \longrightarrow0.
\label{eq:PSnonlinear}
\end{align}
All exponents in this estimate are finite, and those to which a weighted embedding is applied are at least $n$.

Testing $J_\lambda'(u_k)$ with $u_k-u$ and using \eqref{eq:PSnonlinear} gives
\[
 \int_{\R^n}\bigl(|\nabla u_k|^{n-2}\nabla u_k\cdot\nabla(u_k-u)
 +V|u_k|^{n-2}u_k(u_k-u)\bigr)\dd\longrightarrow0.
\]
The linear functional
\[
 h\longmapsto\int_{\R^n}
 \bigl(|\nabla u|^{n-2}\nabla u\cdot\nabla h
       +V|u|^{n-2}uh\bigr)\dd
\]
is bounded on $E$ by H\"older's inequality, with absolute value at most $\|u\|^{n-1}\|h\|$. Its
value at $h=u_k-u$ therefore tends to zero by weak convergence. Subtracting this expression from the
preceding limit and applying the monotonicity inequality, valid for $n\geq2$,
\[
 (|\xi|^{n-2}\xi-|\eta|^{n-2}\eta)\cdot(\xi-\eta)
 \geq2^{2-n}|\xi-\eta|^n,
\]
to both the gradients and the scalar values of the functions yields
\begin{align*}
 2^{2-n}\|u_k-u\|^n
 &\leq\int_{\R^n}
 \bigl(|\nabla u_k|^{n-2}\nabla u_k-|\nabla u|^{n-2}\nabla u\bigr)
 \cdot(\nabla u_k-\nabla u)\dd\\
 &\quad+\int_{\R^n}V
 \bigl(|u_k|^{n-2}u_k-|u|^{n-2}u\bigr)(u_k-u)\dd
 \longrightarrow0.
\end{align*}
Since $V\geq0$, multiplication by the potential preserves the scalar inequality, including where the
potential vanishes. Hence $u_k\to u$ strongly in $E$, and Lemma~\ref{lem:functional} yields
$J_\lambda(u)=c$ and $J_\lambda'(u)=0$, completing the proof.
\end{proof}

We can now prove the existence result and its quantitative consequences.

\begin{proof}[Proof of Theorem~\ref{thm:equation}]
By Lemma~\ref{lem:geometry}, the functional exhibits mountain-pass geometry and $c_\lambda\to0$.
Choose $\lambda_0$ such that $0<c_\lambda<c_*$ for every $\lambda>\lambda_0$. The mountain-pass
theorem of Ambrosetti and Rabinowitz~\cite{AmbrosettiRabinowitz} yields a Palais--Smale sequence at
level $c_\lambda$. Lemma~\ref{lem:PS} provides a strongly convergent subsequence whose limit
$u_\lambda$ satisfies
\[
 J_\lambda'(u_\lambda)=0,\qquad J_\lambda(u_\lambda)=c_\lambda>0.
\]
Therefore $u_\lambda\neq0$ and \eqref{eq:weaksolution} holds. Testing with
$u_\lambda^-:=\max\{-u_\lambda,0\}$ gives $-\|u_\lambda^-\|^n=0$, since $f$ vanishes on
$(-\infty,0]$. Thus $u_\lambda\geq0$.
\end{proof}

For the large-parameter existence conclusion alone, any uniform singular exponential bound with
positive coefficient on the unit energy ball suffices. Indeed, if $\|u\|_{W^{1,n}}\leq C\|u\|$ with
the standard full Sobolev norm, the constant-potential inequality of Adimurthi--Yang~\cite{AdimurthiYang} gives such a bound with coefficient $C^{-n'}\alpha_{n,\beta}$.
Lemmas~\ref{lem:geometry} and~\ref{lem:PS} then give compactness at the mountain-pass level for
large $\lambda$ with the corresponding smaller threshold when $C>1$. Using Theorem~\ref{thm:sharp}
yields the sufficient range \eqref{eq:PSthreshold}.

\begin{proof}[Proof of Corollary~\ref{cor:model}]
Set $f(t)=t^{p-1}e^{\alpha_0t^{n'}}$ for $t>0$ and $f(t)=0$ for $t\leq0$. Conditions \eqref{ass:f1}
and \eqref{ass:f2} follow immediately from $p>n$. Since the exponential factor is increasing,
\[
 F(t)=\int_0^t s^{p-1}e^{\alpha_0s^{n'}}\,ds
 \leq\frac{t^p}{p}e^{\alpha_0t^{n'}}=\frac{tf(t)}p,
 \qquad t>0.
\]
Thus \eqref{ass:f3} holds with $\mu=p$. Also $F(t)\geq t^p/p$, so
\[
 J_\lambda(t\varphi)\leq\frac{t^n}{n}\|\varphi\|^n
 -\frac{\lambda t^p}{p}\int_{\R^n}|x|^{-\beta}\varphi^p\dd.
\]
The maximum of the right-hand side occurs at
\[
 t^{p-n}=\frac{\|\varphi\|^n}
 {\lambda\displaystyle\int_{\R^n}|x|^{-\beta}\varphi^p\dd}.
\]
Substituting this value into \eqref{eq:raylevel} proves \eqref{eq:modelbound}. The bound tends to
zero, so Theorem~\ref{thm:equation} applies for sufficiently large $\lambda$.
\end{proof}

\begin{remark}\label{rem:smallsolutions}
For the mountain-pass solutions in Theorem~\ref{thm:equation}, condition~\eqref{ass:f3} gives the
quantitative bound
\begin{equation}
 \left(\frac1n-\frac1\mu\right)\|u_\lambda\|^n
 \leq J_\lambda(u_\lambda)-\frac1\mu J_\lambda'(u_\lambda)u_\lambda
 =c_\lambda\longrightarrow0.
\label{eq:solutionsmallnorm}
\end{equation}
Hence $u_\lambda\to0$ strongly in $E$ as $\lambda\to\infty$, and therefore in the unweighted and
weighted $L^q$ spaces for every finite $q\geq n$.

For the model in Corollary~\ref{cor:model}, there are constants $C_1,C_2>0$ such that, for all
sufficiently large $\lambda$,
\begin{equation}
 C_1\lambda^{-1/(p-n)}\leq\|u_\lambda\|
 \leq C_2\lambda^{-1/(p-n)}.
\label{eq:modelnormrate}
\end{equation}
The upper bound follows from \eqref{eq:modelbound} and \eqref{eq:solutionsmallnorm}, with $\mu=p$.
For the lower bound, fix a sufficiently small $\rho>0$. Splitting the exponential into its first
$n-1$ terms and the remainder $\Phi_n$ and estimating as in \eqref{eq:smallnormestimate} gives
\begin{equation}
 \int_{\R^n}|x|^{-\beta}|u|^p e^{\alpha_0|u|^{n'}}\dd
 \leq C\|u\|^p\quad\text{if }\|u\|\leq\rho.
\label{eq:modelsmallgrowth}
\end{equation}
Indeed, the polynomial terms are controlled by the weighted embeddings, while H\"older's inequality
and the endpoint estimate control the remainder when $2\alpha_0\rho^{n'}<\alpha_{n,\beta}$. Since
\eqref{eq:solutionsmallnorm} ensures $\|u_\lambda\|\leq\rho$ for large $\lambda$, testing the
equation with $u_\lambda$ yields
\[
 \|u_\lambda\|^n
 =\lambda\int_{\R^n}|x|^{-\beta}u_\lambda^p
 e^{\alpha_0u_\lambda^{n'}}\dd
 \leq C\lambda\|u_\lambda\|^p.
\]
Since $u_\lambda\neq0$, division by $\|u_\lambda\|^n$ proves the lower bound in
\eqref{eq:modelnormrate}. These estimates apply to the mountain-pass solutions constructed above,
without any uniqueness assumption.
\end{remark}

\begin{remark}\label{rem:zerogeometry}
The potential in \eqref{eq:decreasingexample} satisfies the hypotheses for both the critical
extremal and the large-parameter solution. For the equation the zero set may also be disconnected.
Given any nonempty compact set $K\subset\R^n$, the potential
\begin{equation}
 V(x)=\frac{\operatorname{dist}(x,K)^2}{1+\operatorname{dist}(x,K)^2}
 \label{eq:compactzeroset}
\end{equation}
is continuous, has zero set $K$ and tends to one at infinity. Thus \eqref{ass:V} holds. The sharp
bound, subcritical attainment and the mountain-pass solution apply whether or not $0\in K$.
\end{remark}

\noindent\textbf{Acknowledgements.} This work was supported by the National Natural Science Foundation of China (Nos. 12201089, 12371051, and 12526101), the Chongqing Normal University Foundation (No. BWQB2023006), and the Open Research Fund of the Key Laboratory of Nonlinear Analysis and Applications, Ministry of Education, Central China Normal University (No. NAA2025ORG006).



\end{document}